\documentclass[12pt]{article}
\usepackage{amsmath,amsthm,amssymb}
\begin{document}
\newtheorem{thm}{Theorem}
\newtheorem{pro}[thm]{Proposition}
\newtheorem{cor}[thm]{Corollary}
\newtheorem{lem}[thm]{Lemma}
\newtheorem{df}[thm]{Definition}
\newtheorem{rem}[thm]{Remark}
\newtheorem{conj}[thm]{Conjecture}
\newcommand{\End}{\mathop{\mathrm{End}}\nolimits}
\newcommand{\Id}{\mathop{\mathrm{Id}}\nolimits}
\newcommand{\Ric}{\mathop{\mathrm{Ric}}\nolimits}
\title{Uniqueness of Hamiltonian volume minimizing Lagrangian submanifolds
which are Hamiltonian isotopic to ${\mathbb R}P^{n}$ in ${\mathbb C}P^{n}$
}
\author{Hiroshi Iriyeh}
\date{}
\maketitle

\section{Introduction}

\indent\indent
The equator on $S^2$ has the least length among all its images under
area bisecting deformations.
A symplectic higher dimensional generalization of this result was obtained
by Y.-G.\ Oh \cite{Oh} and B.\ Kleiner in 1990.
\begin{thm}[Kleiner-Oh]
The standard real form ${\mathbb R}P^{n} \subset {\mathbb C}P^{n}$
has the least volume among all its images under Hamiltonian isotopies.
\end{thm}
In other words, for any Hamiltonian diffeomorphism
$\phi \in \mathrm{Ham}({\mathbb C}P^{n})$ of the $n$-dimensional complex
projective space ${\mathbb C}P^{n}$, we have
 \begin{eqnarray*}
  \mathrm{vol}(\phi({\mathbb R}P^{n})) \geq \mathrm{vol}({\mathbb R}P^{n}).
 \end{eqnarray*}
A minimal (or Hamiltonian-minimal) Lagrangian submanifold with such a property
is said to be {\it Hamiltonian volume minimizing}.

The purpose of this short note is to prove the uniqueness of Hamiltonian
volume minimizing Lagrangian submanifolds which are Hamiltonian isotopic to
${\mathbb R}P^{n}$ in ${\mathbb C}P^{n}$ modulo isometric group actions.
More precisely,
\begin{thm}
Let $\phi$ be a Hamiltonian diffeomorphism of ${\mathbb C}P^{n}$ with the
standard Fubini-Study K\"{a}hler form $\omega$.
If $\mathrm{vol}(\phi({\mathbb R}P^{n})) = \mathrm{vol}({\mathbb R}P^{n})$
holds, then $\phi$ must be an isometry of ${\mathbb C}P^{n}$.
\end{thm}
\begin{rem}
In the case of a totally geodesic Lagrangian torus in $S^2(1) \times S^2(1)$,
an analogous result of the above theorems was proved in {\rm \cite{IOS}} and
{\rm \cite{IOS2}}.
\end{rem}
The author is grateful to Professor T.\ Otofuji for helpful comments.

\section{Review on Kleiner-Oh's method}

\indent\indent
In this section, we review a proof of theorem 1 by Kleiner and Oh.
This is a combination of the Lagrangian intersection theorem in
symplectic geometry and Poincar\'e formula in integral geometry.
Let $({\mathbb C}P^{n},\omega,J)$ be the $n$-dimensional complex projective
space with the Fubini-Study K\"{a}hler form $\omega$ and the standard
complex structure $J$.
\begin{thm}[Givental \cite{Gi}]
For any Hamiltonian diffeomorphism $\phi$ of ${\mathbb C}P^{n}$ such that
${\mathbb R}P^{n}$ and $\phi({\mathbb R}P^{n})$ intersect transversely,
the inequality
 \begin{eqnarray}
  \sharp ({\mathbb R}P^{n} \cap \phi({\mathbb R}P^{n})) \geq n+1
 \end{eqnarray}
holds.
\end{thm}
Next, we explain the Poincar\'e formula for Lagrangian submanifolds in
${\mathbb C}P^{n}$.
\begin{thm}[Kleiner/Howard \cite{Howard}]
Let ${\mathbb C}P^{n}=U(n+1)/U(1) \times U(n)$ be the complex projective
space.
Let $L_{1}$ and $L_{2}$ be compact Lagrangian submanifolds in
${\mathbb C}P^{n}$.
Then we have
 \begin{eqnarray}
  \int_{U(n+1)} \sharp(L_{1} \cap gL_{2})d\mu(g)=
  C \mathrm{vol}(L_{1})\mathrm{vol}(L_{2}) 
 \end{eqnarray}
 where
  \begin{eqnarray*}
   C=\frac{(n+1)\mathrm{vol}(U(n+1))}{\mathrm{vol}({\mathbb R}P^{n})^2}.
  \end{eqnarray*}
\end{thm}
These two theorems yield the Hamiltonian volume minimizing property of
${\mathbb R}P^{n}$ as follows.
By equation $(2)$, we have
\begin{eqnarray}
 \frac{(n+1)\mathrm{vol}(U(n+1))}{\mathrm{vol}({\mathbb R}P^{n})^2}
  \mathrm{vol}({\mathbb R}P^{n})\mathrm{vol}(\phi({\mathbb R}P^{n}))\nonumber\\
 =\int_{U(n+1)} \sharp({\mathbb R}P^{n} \cap g\phi({\mathbb R}P^{n}))d\mu(g).
\end{eqnarray}
Since $g \in U(n+1)$ is a Hamiltonian diffeomorphism,
by the group structure of $\mathrm{Ham}({\mathbb C}P^{n})$,
$g\phi$ is also a Hamiltonian diffeomorphism.
By equation $(1)$, we have
\begin{eqnarray}
 {\rm RHS\ of\ } (3) \geq \int_{U(n+1)} (n+1)d\mu(g)
 =(n+1)\mathrm{vol}(U(n+1)).
\end{eqnarray}
Hence,
\begin{eqnarray*}
 \mathrm{vol}(\phi({\mathbb R}P^{n})) \geq \mathrm{vol}({\mathbb R}P^{n}).
\end{eqnarray*}

\section{Proof}

\indent\indent
First of all, we exhibit the key fact to prove theorem 2.
\begin{thm}[Oh \cite{Oh2}]
Let $L$ be an embedded closed Lagrangian submanifold in ${\mathbb C}P^{n}$.
If $L$ satisfies that
 \begin{eqnarray*}
  \sharp(L \cap gL)=SB(L,{\mathbb Z}_2)
 \end{eqnarray*}
for all $g \in U(n+1)$ such that $L$ and $gL$ intersect transversely,
then $L$ is the standard totally geodesic ${\mathbb R}P^{n}$ in
${\mathbb C}P^{n}$.
Here, $SB(L,{\mathbb Z}_2)$ denotes the sum of ${\mathbb Z}_2$-Betti
numbers of $L$.
\end{thm}

Assume that
$\mathrm{vol}(\phi({\mathbb R}P^{n})) = \mathrm{vol}({\mathbb R}P^{n})$
for $\phi \in \mathrm{Ham}({\mathbb C}P^{n})$.
Then we have
\begin{eqnarray}
 \int_{U(n+1)} \sharp({\mathbb R}P^{n} \cap g\phi({\mathbb R}P^{n}))d\mu(g)
  =\int_{U(n+1)} (n+1)d\mu(g).
\end{eqnarray}
On the other hand, by assumption,
\begin{eqnarray}
 \lefteqn{\int_{U(n+1)}
  \sharp(\phi({\mathbb R}P^{n}) \cap g\phi({\mathbb R}P^{n}))d\mu(g)}
  \nonumber\\
 &=& C \mathrm{vol}(\phi({\mathbb R}P^{n}))
      \mathrm{vol}(\phi({\mathbb R}P^{n}))\nonumber\\
 &=& C \mathrm{vol}({\mathbb R}P^{n})
      \mathrm{vol}(\phi({\mathbb R}P^{n}))\nonumber\\
 &=& \int_{U(n+1)}
      \sharp({\mathbb R}P^{n} \cap g\phi({\mathbb R}P^{n}))d\mu(g).
\end{eqnarray}
By equations (5) and (6), we have
\begin{eqnarray*}
 \int_{U(n+1)}
  \{ \sharp(\phi({\mathbb R}P^{n}) \cap g\phi({\mathbb R}P^{n}))-(n+1) \}
   d\mu(g)=0.
\end{eqnarray*}
By theorem 4, for almost every $g \in U(n+1)$, we have
\begin{eqnarray}
 \sharp(\phi({\mathbb R}P^{n}) \cap g\phi({\mathbb R}P^{n}))=n+1.
\end{eqnarray}

Take $g \in U(n+1)$ which satisfies that $\phi({\mathbb R}P^{n})$ and
$g\phi({\mathbb R}P^{n})$ intersect transversely.
If $\sharp(\phi({\mathbb R}P^{n}) \cap g\phi({\mathbb R}P^{n})) > n+1$,
then this inequality holds for all $g' \in U(n+1)$ near $g$.
This contradicts equation (7).
Hence, $\sharp(\phi({\mathbb R}P^{n}) \cap g\phi({\mathbb R}P^{n}))=n+1$
for all $g \in U(n+1)$ such that $\phi({\mathbb R}P^{n})$ and
$g\phi({\mathbb R}P^{n})$ intersect transversely.

By theorem 6, the Lagrangian submanifold $\phi({\mathbb R}P^{n})$
must be the standard totally geodesic
${\mathbb R}P^{n} \subset {\mathbb C}P^{n}$.
Therefore, $\phi \in U(n+1)$. \hfill \qed

\begin{flushleft}
{\sc School of Engineering\\
Tokyo Denki University\\
Kanda-Nishiki-Cho, Chiyoda-Ku\\
Tokyo, 101-8457\\
Japan}

{\it e-mail} : {\tt hirie@im.dendai.ac.jp}
\end{flushleft}
\end{document}